\documentclass[12pt,leqno]{article}
\usepackage{amsmath,amssymb,amsfonts,amsthm,latexsym,amstext,amscd}
\usepackage[english]{babel}
\usepackage{fancyhdr}

\newtheorem{theorem}{Theorem}

\newtheorem{lemma}[theorem]{Lemma}

\newtheorem*{prop}{Proposition}
\newtheorem*{theo}{Theorem}

\theoremstyle{definition}

\newtheorem*{defi}{Definition}
\newtheorem*{rema}{Remark}

\newcommand{\beq}{\begin{equation}}
\newcommand{\eeq}{\end{equation}}

\begin{document}

\title
{On the ratio of consecutive gaps between primes}
\author
{J\'anos Pintz\thanks{Supported by OTKA Grants NK104183, K100291 and ERC-AdG.~321104.}\\
\\
\emph{Dedicated to Helmut Maier on the occasion of his 60\textsuperscript{th} birthday}}

\date{}

\numberwithin{equation}{section}


\maketitle

\section{Introduction}
\label{sec:1}

The difference between the consecutive primes, the expression
\beq
\label{eq:1.1}
d_n = p_{n + 1} - p_n,
\eeq
where $\mathcal P = \{p_i\}_{i = 1}^\infty$ denotes the set of primes, has been investigated probably since the time of the Greeks.
The Twin Prime Conjecture asserts
\beq
\label{eq:1.2}
d_n = 2 \ \text{ infinitely often.}
\eeq

This conjecture settles at the conjectural level the small values of $d_n$.
Concerning the large values even the suitable conjecture is not completely clear.
However, it seems to be that Cram\'er's conjecture \cite{Cra1, Cra2}
\beq
\label{eq:1.3}
\limsup\limits_{n \to \infty} \frac{d_n}{\log^2 n} = C_0 = 1
\eeq
is near to the truth.
Granville suggested \cite{Gra1, Gra2}, just based on the famous matrix method of Helmut Maier that the correct value of $C$ is instead of $1$ slightly larger
\beq
\label{eq:1.4}
C \geqslant 2e^{-\gamma} > 1.
\eeq
However, most mathematicians agree that the correct maximal order of $d_n$ should be $(\log n)^{2 + o(1)}$.

In the present work we give a short overview of the older and the recent new breakthrough results of small gaps between primes, of the present state of art of long differences.
Our final result will be a common generalization of them which deals with a problem raised first by Erd\H{o}s and Tur\'an about small and large values of
\beq
\label{eq:1.5}
\frac{d_{n + 1}}{d_n}.
\eeq

As it will be clear from the history and from the proof of our new result we will often refer to important and ingenious results of Helmut Maier, which often (but not exclusively) use his famous matrix method.

Since there were more than 40 papers discussing small and large values of $d_n$ (or their relations) we will be relatively brief and try to concentrate on the most important developments.

Due to the prime number theorem,
\beq
\label{eq:1.6}
\frac1{N} \sum_{n = 1}^N d_n = \frac{p_{N + 1} - 2}{N} \sim \log N
\eeq
we have trivially
\beq
\label{eq:1.7}
\Delta := \liminf_{n \to \infty} \frac{d_n}{\log n} \leqslant 1 \leqslant \limsup_{n \to \infty}  \frac{d_n}{\log n} =: \lambda.
\eeq

\subsection*{A) Large gaps between primes}

The first non-trivial result, $\lambda > 1$, was reached in 1929, when Backlund \cite{Bac} proved $\lambda \geqslant 2$.
One year later Brauer and Zeitz \cite{BZ} improved this to $\lambda \geqslant 4$.
Another year later E. Westzynthius \cite{Wes} showed already $\lambda = \infty$, more precisely
\beq
\label{eq:1.8}
\limsup_{n \to \infty} \frac{d_n / \log n}{\log_3 n / \log_4 n} \geqslant 2e^\gamma,
\eeq
where $\log_\nu n$ denotes the $\nu$-fold iterated logarithmic function.

Erd\H{o}s \cite{Erd1} improved this in 1935 to
\beq
\label{eq:1.9}
\limsup_{n \to \infty} \frac{d_n/\log n}{\log_2 n / (\log_3 n)^2} > 0.
\eeq

Finally, three years later Rankin \cite{Ran1} succeeded in showing the function which is even currently, after more than 75 years the best result, apart from the constant:
\beq
\label{eq:1.10}
\limsup_{n \to \infty} \frac{d_n / \log n}{\log_2 n \log_4 n / (\log_3 n)^2} \geqslant \frac13.
\eeq

Since, apart from the improvements of the constant $1/3$ to $e^\gamma$ (in works of Ricci, Rankin and Sch\"onhage between 1952--1963), the inequality \eqref{eq:1.10} remained the strongest, Erd\H{o}s offered in 1979, at a conference in Durham a prize of USD 10,000 for a proof of \eqref{eq:1.10} with an arbitrarily large constant $C$ in place of~$1/3$; the highest prize ever offered by him for a mathematical problem.
However, this did not help either in the past 35 years.
Nevertheless, strong analytic methods introduced by Helmut Maier and Carl Pomerance, combined with the original Erd\H{o}s--Rankin sieve procedure helped to prove \cite{MaP} in 1990 \eqref{eq:1.10} with $C = 1.3126\dots e^\gamma$.

Finally, the current best result, the relation \eqref{eq:1.10} with $1/3$ replaced by
\beq
\label{eq:1.11}
C = 2e^\gamma
\eeq
was reached by J. Pintz \cite{Pin1}.
The methods to reach it involve besides the classical sieve methods and the large sieve used by Helmut Maier and Carl Pomerance, a pure graph-theoretical result, which, however, is proved in \cite{Pin1} using a so-called semi-random method of E.\ Szemer\'edi.

\subsection*{B) Chains of Large Gaps Between Consecutive Primes}

In 1949 Erd\H{o}s \cite{Erd4} proposed the problem whether $k$ consecutive prime gaps $d_{n + 1}, \dots, d_{n + k}$ can
be simultaneously much larger than its mean value, i.e.\ whether
\beq
\label{eq:1.12}
\limsup_{n \to \infty} \frac{\min (d_{n + 1}, \dots, d_{n + k})}{\log n} = \infty,
\eeq
and succeeded in showing this for $k = 2$ in the same work.

30 years later Helmut Maier \cite{Mai1} introduced his famous matrix method and showed this conjecture even in a stronger form when $\log n$ is replaced by the Erd\H{o}s--Rankin function (cf.\ \eqref{eq:1.10}).
He proved for any $k$
\beq
\label{eq:1.13}
\limsup_{n \to \infty} \frac{\min(d_{n + 1}, \dots, d_{n + k})/\log n}{\log_2 n \log_4 n / (\log_3 n)^2} > 0.
\eeq

\subsection*{C) Small Gaps Between Consecutive Primes}

Unlike the quick progress in the case of large gaps an analogue of the early result $\lambda \geqslant 4$ of Brauer--Zeitz from 1930, i.e.\ $\Delta \leqslant 1/4$ was essentially the best result still even 75 years later before 2005, and it was reached by Helmut Maier \cite{Mai2} also by his matrix method in 1985.

The first, although unpublished and conditional result was reached by Hardy and Littlewood (see \cite{Ran2}) in 1926 who showed that the Generalized Riemann Hypothesis (GRH) implies $\Delta \leqslant 2/3$.
However, the first non-trivial unconditional result was proved by Erd\H{o}s \cite{Erd2},
\beq
\label{eq:1.14}
\Delta \leqslant 1 - c_1,
\eeq
with an unspecified explicitly calculable constant $c_1 > 0$.

After much work, calculating and improving $c_1$, the next breakthrough came by the large sieve of Bombieri \cite{Bom} and Vinogradov \cite{Vin} which enabled Bombieri and Davenport \cite{BD} to eliminate GRH and (incorporating also Erd\H{o}s' ideas into their work) to show unconditionally
\beq
\label{eq:1.15}
\Delta \leqslant (2 + \sqrt{3}) / 8 = 0.466\dots\,.
\eeq

After five further improvements of Piltjai, Huxley and Fouvry--Grupp this was reduced to 0.4342.
The next big step was the mentioned result of Helmut Maier \cite{Mai2}, the inequality
\beq
\label{eq:1.16}
\Delta \leqslant 0.2486.
\eeq

Twenty years later D. Goldston, J. Pintz and C. Y{\i}ld{\i}r{\i}m succeeded in reaching the optimal value
\beq
\label{eq:1.17}
\Delta = 0
\eeq
(see Primes in Tuples I and III in \cite{GPY1} and \cite{GPY2}).
Further they showed that $d_n$ can be much smaller than $\log n$, in fact it was proved in \cite{GPY3} that
\beq
\label{eq:1.18}
\liminf \frac{d_n}{(\log n)^c} = 0 \ \text{ for any } \ c > 1/2.
\eeq

This was further improved by J. Pintz \cite{Pin4} for any $c > 3/7$ which in some sense was the limit of the original GPY method as proved by B. Farkas, J. Pintz and Sz. Gy. R\'ev\'esz \cite{FPR}.

On the other hand, already the first work \cite{GPY1}, describing $\Delta = 0$ in details contained a conditional theorem.
In order to formulate it we will describe the notion of admissibility for a $k$-tuple $\mathcal H = \{h_i\}_{i = 1}^k$ of
distinct integers.

\begin{defi}
$\mathcal H = \{h_i\}_{i = 1}^k$ ($h_i \neq h_j$ for $i \neq j$) is called admissible if for any prime $p$ the set $\mathcal H$ does not occupy all residue classes $\text{\rm mod }p$.

This is equivalent with the formulation that $\prod\limits_{i = 1}^k (n + h_i)$ has no fixed prime divisor.
\end{defi}

\begin{rema}
Although we usually suppose that $h_i \geqslant 0$ this is clearly not necessary in view of the consequence which is translation invariant.
\end{rema}

We also introduced two connected conjectures which were named after Dickson, Hardy and Littlewood and which represent actually a weaker form of Dickson's prime $k$-tuple conjecture \cite{Dic}.

\bigskip
\noindent
{\bf Conjecture DHL$(k, k_0)$.}
{\it If $\mathcal H_k$ is admissible of size $k$, then
the translated sets $n + \mathcal H_k$ contain for infinitely many $n$ values at least $k_0$ primes if $k > C(k_0)$.}

\bigskip
The special case $k_0 = 2$ had special attention because of his connection with the

\bigskip
\noindent
{\bf Bounded Gap Conjecture.}
{\it $\liminf\limits_{n \to \infty} d_n \leqslant C$ with a suitable absolute constant $C$.}

\bigskip
The mentioned connection is the simple

\begin{prop}
If Conjecture {\rm DHL}$(k,2)$ is true for some $k$, then the Bounded Gap Conjecture is true.
\end{prop}

To formulate our conditional result we still need the following

\begin{defi}
A number $\vartheta$ is called a level of distribution of primes if for any $A, \varepsilon > 0$ we have
\beq
\label{eq:1.19}
\sum_{q \leqslant X^{\vartheta - \varepsilon}} \max_{\substack{a\\
(a,q) = 1}}
\biggl| \sum_{\substack{p \leqslant X\\
p \equiv a (\text{\rm mod }q)}} \log p - \frac{X}{\varphi(q)} \biggr| \leqslant \frac{C(A, \varepsilon)X}{(\log X)^A}.
\eeq
\end{defi}

The result that \eqref{eq:1.19} holds with $\vartheta = 1/2$ is the famous Bombieri--Vi\-no\-g\-ra\-dov theorem (\cite{Bom}, \cite{Vin}), while Elliott--Halberstam \cite{EH} conjectured that even $\vartheta = 1$ is a level of distribution.

We can introduce for any $\vartheta \in (1/2, 1]$ the

\bigskip
\noindent
{\bf Conjecture EH$(\vartheta)$.}
{\it \eqref{eq:1.19} is true for $\vartheta$, i.e.\ $\vartheta$ is a level of distribution of the primes.}

\bigskip
We showed in our original work

\begin{theo}[{\cite{GPY1}}]
If {\rm EH}$(\vartheta)$ is true for any $\vartheta > 1/2$, then {\rm DHL}$(k, 2)$ is true for $k > C_1(\vartheta)$ and consequently,
$\liminf d_n \leqslant C_2(\vartheta)$.
\end{theo}

Soon after this, Y. Motohashi and J. Pintz \cite{MP}, (MR 2414788 (2009d:1132), arXiv: math/0602599, Feb 27, 2006) showed in the work entitled ``A Smoothed GPY sieve'' that \eqref{eq:1.19} can be substituted with the weaker condition that it holds for smooth moduli of $q$ in \eqref{eq:1.19} (by that we mean that for moduli $q$ having all their prime factors below $q^b$ -- or even $X^b$ -- with an arbitrarily fixed $b$)
and the maximum taken over all residue classes $a$ with $(a, q) = 1$ can be reduced to those satisfying $\prod\limits_{i = 1}^k (a + h_i ) \equiv 0 \ (\text{\rm mod }q)$ which are trivially the only cases appearing in the proof.

Finally, Y. T. Zhang \cite{Zha} showed that for $b = \frac1{292}$ the above condition holds with $\vartheta = \frac12 + \frac1{584}$.

This led to the very recent result.

\begin{theo}[{\cite{Zha}}]
{\rm DHL}$(3.5 \cdot 10^6, 2)$ is true and consequently $\liminf\limits_{n \to \infty} d_n \leqslant 7 \cdot 10^7$.
\end{theo}

\begin{rema}
Y. Zhang attributes the result \cite{MP} to himself (and proves it again in his work \cite{Zha}) despite the fact that the authors Motohashi and Pintz called his attention in four subsequent e-mails to their work
and asked him to mention this fact in his manuscript
 in May--June 2013, when his manuscript appeared first electronically on the webpage of Annals of Mathematics.
 He completely ignored and left unanswered three of them and answered the fourth one in one line, refusing to add anything on it, based on his assertion that ``when preparing my manuscript I had not read your paper.''
 Finally the printed version of his work appeared without any reference to \cite{MP} (see \cite{Zha}),
 despite the fact that he described in an interview that he got the idea to use a smoothed GPY sieve on July 3rd, 2012 and that was crucial to his solution \cite{Kla}.
 This was more than 6 years later than the appearance of ``A smoothed GPY sieve'' on arXiv on February 27, 2006.
\end{rema}

Some months later the joint effort of many mathematicians showed a stronger form of Zhang's theorem in the Polymath 8A project led by T. Tao.

\begin{theo}[Polymath 8A]
{\rm DHL}$(632, 2)$ is true and consequently $\liminf\limits_{n \to \infty} d_n \leqslant 4680$.
\end{theo}

Another improved version, proved independently by another refinement of the GPY method (in spirit closer to an elementary version of the method worked out in collaboration with S. W. Graham, see \cite{GGPY} and which was also close to the first attempt of Goldston and Y{\i}ld{\i}r{\i}m \cite{GY} which finally led to $\Delta \leqslant 1/4$ in their version) was reached by J. Maynard \cite{May}, which showed the even stronger

\begin{theo}[{\cite{May}}]
{\rm DHL}$(105, 2)$ is true and consequently $\liminf\limits_{n \to \infty} d_n \leqslant 600$.
\end{theo}

T. Tao used the same approach independently and simultaneously with Maynard and with help of his Polymath Project 8B (involving further new theoretical ideas and a huge number of computations) they showed

\begin{theo}[Polymath~8B]
$\text{\rm DHL}(50, 2)$ is true and consequently $\liminf\limits_{n \to \infty} d_n \leqslant 246$.
\end{theo}

It is interesting to note that the above methods did not need the results or ideas of Motohashi--Pintz and Zhang, neither any weaker form of a result of type $\vartheta > 1/2$.

In complete contrast to this, the Maynard--Tao method shows the existence of infinitely many bounded gaps (although with weaker numerical bounds than 600 or 246) with any fixed positive distribution level of the primes.
The first result of this type was proved by A. R\'enyi in 1947--48 \cite{Ren}.

\subsection*{D) Chains of Bounded Gaps Between Consecutive Primes}

The original GPY method \cite{GPY1} furnished only under the very strong original Elliott--Halberstam Conjecture EH(1) the bound
\beq
\label{eq:1.20}
\Delta_2  := \liminf\limits_{n \to \infty} \frac{p_{n + 2} - p_n}{\log p_n} = 0.
\eeq
Thus it failed to show on EH even DHL$(k, 3)$ for some~$k$.
The additional ideas of Motohashi--Pintz and Zhang helped neither.

Already Erd\H{o}s mentioned \cite{Erd4} as a conjecture
\beq
\label{eq:1.21}
\liminf\limits_{n \to \infty} \frac{\min(d_n, d_{n + 1})}{\log n} < 1
\eeq
which was proved in the mentioned work of Helmut Maier \cite{Mai2}.
In the same work he has shown
\beq
\label{eq:1.22}
\Delta_r = \liminf\limits_{n \to \infty} \frac{p_{n + r} - p_n}{\log n} \leqslant e^{-\gamma} \left(r - \frac58 + o(1)\right) \quad (r \to \infty).
\eeq
This was improved in \cite{GPY2} to
\beq
\label{eq:1.23}
\Delta_r \leqslant e^{-\gamma} \left(\sqrt{r} - 1\right)^2.
\eeq
However, as mentioned already, even \eqref{eq:1.20} was open unconditionally.
In this aspect the Maynard--Tao method was much more successful, yielding

\begin{theo}[Maynard--Tao]
$\liminf\limits_{n \to \infty} (p_{n + r} - p_n) \leqslant Ce^{4r}$ for any $r$ with an absolute constant~$C$.
\end{theo}

\section{Consecutive values of $d_n$.\
Problems of Erd\H{o}s, Tur\'an and P\'olya}
\label{sec:2}

The questions discussed in Section~\ref{sec:1} referred to single or consecutive small values of $d_n$ or to analogous  problems dealing with solely large values of~$d_n$.
In 1948 Erd\H{o}s and Tur\'an \cite{ET} showed that $d_{n + 1} - d_n$ changes sign infinitely often.
After this, still in the same year, Erd\H{o}s \cite{Erd3} proved that
\beq
\label{eq:2.1}
\liminf \frac{d_{n + 1}}{d_n} < 1 < \limsup \frac{d_{n + 1}}{d_n}.
\eeq
He mentioned 60 years ago \cite{Erd5}:
``One would of course conjecture that
\beq
\label{eq:2.2}
\liminf\limits_{n \to \infty} \frac{d_{n + 1}}{d_n} = 0 \ \text{ and } \ \limsup_{n \to \infty} \frac{d_{n + 1}}{d_n} = \infty,
\eeq
but these conjectures seem very difficult to prove.''
Based on a generalization of the method of Zhang \cite{Zha} the author proved \eqref{eq:2.2} in \cite{Pin2}.

In the mentioned work of Erd\H{o}s and Tur\'an \cite{ET} they also asked for a necessary and sufficient condition that
\beq
\label{eq:2.3}
\sum_{i = 1}^k a_i p_{n + i}
\eeq
should change sign infinitely often as $n \to \infty$.
They observed that the condition
\beq
\label{eq:2.4}
\sum_{i = 1}^k a_i = 0
\eeq
is clearly necessary.
Using \eqref{eq:2.4} one can reformulate the problem and ask for infinitely many sign changes of
\beq
\label{eq:2.5}
\sum_{i = 1}^\ell \alpha_i d_{n + i} = - \sum_{i = 1}^k a_i p_{n + i}
\eeq
if we use the notation
\beq
\label{eq:2.6}
\alpha_j = \sum_{i = 1}^j a_i \quad (j = 1,2, \dots, k - 1), \quad \ell = k - 1.
\eeq

This form shows an observation of P\'olya (see \cite{Erd6}) according to which $\alpha_j$ $(j = 1,\dots, \ell)$ cannot all have the same sign if \eqref{eq:2.3} has infinitely many sign changes.
Erd\H{o}s \cite{Erd6} writes: ``It would be reasonable to conjecture that P\'olya's condition is necessary and sufficient for \eqref{eq:2.5} to change sign infinitely often.
Unfortunately the proof of this is not likely to succeed at the present state of science.''

After this Erd\H{o}s showed \cite{Erd6} that \eqref{eq:2.3}, i.e.\ \eqref{eq:2.5} changes sign infinitely often if
\beq
\label{eq:2.7}
\sum_{i = 1}^\ell \alpha_i = 0, \quad \alpha_\ell \neq 0.
\eeq

The author announced in \cite{Pin3} the proof of the following

\begin{theo}[{\cite{Pin3}}]
The sum \eqref{eq:2.3}, i.e.\ \eqref{eq:2.5}, changes sign infinitely often if at least one of the following conditions holds $(\ell \geqslant 2)$

{\rm (i)} \quad $\biggl|\displaystyle \sum_{i = 1}^\ell \alpha_i \biggr| \leqslant c_0(\ell) \sum_{i = 1}^\ell |\alpha_i|$,

\noindent
with a sufficiently small explicitly calculable constant $c_0(\ell)$ depending on~$\ell$;

\smallskip
{\rm (ii)} if $\exists j \in [1, \ell]$ such that
\beq
\label{eq:2.8}
\sum_{i = 1, i \neq j}^\ell |\alpha_i| < |\alpha_j|, \quad \text{\rm sgn}\, \alpha_i \neq \text{\rm sgn}\, \alpha_j, \ \ i \in [1, \ell] \setminus j;
\eeq

{\rm (iii)} if the Hardy--Littlewood prime $k$-tuple conjecture is true for $k = \ell$.
\end{theo}

Now, (iii) shows that the mentioned conjecture of Erd\H{o}s, namely that P\'olya's trivial necessary condition (i.e.\ that all $\alpha_j$ cannot have the same sign) is probably really a necessary and sufficient condition for \eqref{eq:2.5} to change sign infinitely often.

\section{Results}
\label{sec:3}

In the present work we will show a kind of improvement of the result \eqref{eq:1.10}, which can be considered also as a common generalization of our result \eqref{eq:2.2} and Maynard--Tao's mentioned theorem about the strongest known estimates of chains of small
gaps between consecutive primes.
The result will also give an improvement on the recent proof of the author (\cite{Pin2}) which solved the 60-year-old problem \eqref{eq:2.2} of Erd\H{o}s and improves the result \eqref{eq:1.12} for $k = 2$ due to Erd\H{o}s \cite{Erd4}.

The proof will also show a simple method (already observed in \cite{Pin2} by the author after the proof of Zhang \cite{Zha})
which makes the producing of bounded gaps (or chains of bounded gaps in case of Maynard and Tao) effective, since the original versions used a Bombieri--Vinogradov type theorem which again made use of the ineffective Siegel--Walfisz theorem.

\begin{rema}
The recent work \cite{BFM} shows implicitly a way to make the Maynard--Tao theorem effective.
\end{rema}

Finally we mention that concerning the Theorem \cite{Pin2} of Section~\ref{sec:2} we can further give a very simple sufficient condition for \eqref{eq:2.3}, i.e.\ \eqref{eq:2.5} to change sign infinitely often.
To simplify the problem we can clearly suppose $\alpha_1 \neq 0$, $\alpha_\ell \neq 0$.
In this case the sufficient condition is simply:
\beq
\label{eq:3.1}
\text{\rm sgn }\alpha_1 \neq \text{\rm sgn } \alpha_\ell.
\eeq
(If we do not make the trivial supposition $\alpha_1 \neq 0$, $\alpha_\ell \neq 0$
we can clearly formulate it in the way that the first and last non-zero elements of the sequence $\{\alpha_i\}_{i = 1}^\ell$
should have opposite sign.)

Summarizing, we will prove the following results, using the basic notation $d_n = p_{n + 1} - p_n$ of \eqref{eq:1.1}, $\log_\nu n$ for the $\nu$-fold iterated logarithmic function.

\begin{theorem}
\label{th:1}
Let $k$ be an arbitrary fixed natural number, $c(k)$, $N(k)$ suitable positive, explicitly calculable constants depending only on~$k$.
Then for any $N > N(k)$ there exists an $n \in [N, 2N]$ such that
\beq
\label{eq:3.2}
\frac{d_{n + 1}}{\max(d_n, \dots, d_{n - k + 1})} >
\frac{c(k) \log N \log_2 N \log_4 N}{(\log_3 N)^2}.
\eeq
\end{theorem}

Essentially the same proof gives the analogous result:

\medskip
\noindent
{\bf Theorem 1'.}
{\it Under the conditions of Theorem~\ref{th:1} we have}
\beq
\label{eq:3.3}
\frac{d_{n - k}}{\max(d_n, \dots, d_{n - k + 1})} >
\frac{c(k) \log N \log_2 N \log_4 N}{(\log_3 N)^2}.
\eeq

\begin{theorem}
\label{th:2}
Let $k_0$ be an arbitrary fixed natural number, $c(k_0)$, $N(k_0)$ suitable positive, explicitly calculable constants depending only on~$k_0$.
Then we have a $k \geqslant  k_0$ such that for any $N > N(k_0)$ there exists an $n \in [N, 2N]$ such that
\[
\frac{\min(d_{n - k}, d_{n + 1})}{\max(d_n, \dots, d_{n - k + 1})} > \frac{c(k_0)\log N \log_2 N \log_4 N}{(\log_3 N)^2}.
\]
\end{theorem}

\begin{theorem}
\label{th:3}
Let $\ell \geqslant 2$ be an arbitrary integer, $\alpha_1, \dots, \alpha_\ell$ real numbers with $\alpha_1 \neq 0$, $\alpha_\ell \neq 0$.
Then the expression
\beq
\label{eq:3.4}
\sum_{i = 1}^\ell \alpha_i d_{n + i}
\eeq
changes sign infinitely often as $n \to \infty$ if
\beq
\label{eq:3.5}
\text{\rm sgn }\alpha_1 \neq \text{\rm sgn }\alpha_\ell.
\eeq
\end{theorem}

We remark that Theorem~\ref{th:3} trivially follows from Theorems~\ref{th:1} and~1'.

\section{Proofs}

The proof requires a combination of the methods of Erd\H{o}s--Rankin and that of Maynard--Tao.
Since Theorem~\ref{th:3} follows from Theorems~\ref{th:1} and 1', further Theorem~\ref{th:2} implies Theorems~\ref{th:1} and 1', it is sufficient to prove Theorem~\ref{th:2}.
Concerning Theorem~\ref{th:2} we will show that we will have infinitely many cases when a block of at least $k_0$ consecutive primes in a bounded interval are preceded and followed by two primegaps of length at least
\beq
\label{eq:4.1}
\frac{c_0(k_0) \log N \log_2 N \log_4 N}{(\log_3 N)^2} := c_0(k_0) \log N f(\log N)
\eeq
each, where the corresponding primes $p_\nu'$s satisfy $\nu \in [N, 2N]$.
To be more specific we will choose an arbitrary set of $m$ different primes $\{h_i\}_{i = 1}^m$ with the property that for any $i, j, t \in [1,m]$, $i \neq j$:
\beq
\label{eq:4.2}
m < C_3(m) < h_1 < h_2 < \dots < h_m < C_4(m), \ \ m = \left\lceil C_5 e^{C_6 k_0}\right\rceil, \ \ h_t \nmid h_i - h_j
\eeq
with the absolute constants $C_5, C_6$ and the constants $C_3(m)$ and $C_4(m)$ depending on $m$ to be chosen later.
Denoting
\beq
\label{eq:4.3}
M = \prod_{p_i \leqslant R,\ p_i \neq h_j\ (j = 1,\dots, m)} p_i,
\eeq
we try to determine a residue class $z$ $(\text{\rm mod }M)$ with the property
\beq
\label{eq:4.4}
(z + h_i, M) = 1 \ \ \ (1 \leqslant i \leqslant m),
\eeq
\beq
\label{eq:4.5}
(z \pm \nu, M) > 1 \ \text{ for all } \ \nu \in [0, c_0 (k_0) Rf(R)] \setminus \{h_i\}_{i = 1}^m .
\eeq
We remark that we will have $R f(R) \asymp \log N f(\log N)$.

To describe the procedure briefly we will look for a block of at least $k_0$ primes among the numbers
\beq
\label{eq:4.6}
\ell + h_1, \dots, \ell + h_m, \ \ \ \ell \equiv z \ (\text{\rm mod } M)
\eeq
by the Maynard--Tao method, using the refinements in \cite{BFM}.
This will contain at least $k_0$ consecutive primes in a bounded block, while the two intervals
\beq
\label{eq:4.7}
[\ell - c_0(k_0)Rf(R), \ell - h_m) \ \text{ and } \ (\ell + h_m, \ell + c_0(k_0)Rf(R)]
\eeq
will certainly not contain any prime due to \eqref{eq:4.5}.
Since (cf. \cite{BFM}) we can choose $M$ as large as a suitable small power of $N$, this will prove our theorem as by the Prime Number Theorem we will have
\beq
\label{eq:4.8}
\log N \asymp \log M \sim R, \ \ f(\log N) \sim f(R)
\eeq
and consequently $Rf(R) \asymp \log Nf(\log N)$.

Let us use parameters $0 < v < w < R/2$ to be chosen later, and let
\begin{align}
\label{eq:4.9}
P_1 :&= \prod_{p \leqslant v,\ p \neq h_j \ (1 \leqslant j \leqslant m)} p, \\
\label{eq:4.10}
P_2 :&= \prod_{v < p \leqslant w} p,\\
\label{eq:4.11}
P_3 :&= \prod_{w < p \leqslant R/2} p,\\
\label{eq:4.12}
P_4 :&= \prod_{R/2 < p \leqslant R} p,\\
\label{eq:4.13}
v &= \log^3 R, \ U = c_0(k_0)Rf(R), \ P:= P(R) = P_1 P_2 P_3 P_4,
\end{align}
with a suitably small $c(k_0)$, chosen at the end.
If we have a Siegel-zero, that is, a real primitive character $\chi_1$ $\text{\rm mod }r$, $r \leqslant N$, with a zero $\beta$ of the corresponding $\mathcal L$-function $\mathcal L(s, \chi)$ satisfying
\beq
\label{eq:4.14}
\beta > 1 - \frac{c_7}{\log N},
\eeq
then this character and zero are uniquely determined by the Landau--Page theorem (see e.g.\ \cite[p.\ 95]{Dav}) if $c_7$ is chosen as a suitably small explicitly calculable positive absolute constant.
If $\chi_1$ is real primitive $\text{\rm mod }r$, then $r$ has to be square-free apart from the possibility that the prime $2$ appears on the second or third power in~$r$.
Let us denote the greatest prime factor of~$r$ by~$q$.
Since we have
\beq
\label{eq:4.15}
1 - \beta \gg \frac{\log^2 r}{\sqrt{r}}
\eeq
(see \cite[p.\ 96]{Dav}) where the implied constant is explicitly calculable
we have by \eqref{eq:4.14}
\beq
\label{eq:4.16}
r \gg (\log N)^{3/2}.
\eeq
Consequently, by the Prime Number Theorem and the ``almost square-free'' property of $r$ we have
\beq
\label{eq:4.17}
q \gg \log r \gg \log_2 N \to \infty \ \text{ as } \ N \to \infty.
\eeq
If there is no Siegel-zero, then we will let $q = 1$.
Now we change slightly the definition of $P_i$ and $P$ as to exclude from their defining product the possible divisor~$q$.
Let
\beq
\label{eq:4.18}
P_i^* = P_i /  q \ \text{ if } \ q\mid P_i \ (1 \leqslant i \leqslant 4), \text{ otherwise } P_i^* = P_i, \ P^* = \prod_{i = 1}^4 P_i^*.
\eeq
This change is not necessary for finding $z$ with \eqref{eq:4.4}--\eqref{eq:4.5} but to assure uniform distribution of primes $\equiv z$ $(\text{\rm mod } M)$ in the Maynard--Tao procedure (and to make later the whole procedure effective as well).
This means also that in contrast with the rough description at the beginning of this section the role of $M = P$ will be played actually by
\beq
\label{eq:4.19}
M^* = P^* = P/q \ \text{ if } \ q \mid P\ (\text{otherwise } M^* = M = P = P^*).
\eeq

Returning with this change to our problem of choosing $z \text{ mod } P^*$ let
\beq
\label{eq:4.20}
z \equiv 0 \ (\text{\rm mod }P_1^* P_3^*)
\eeq
while we will choose $z$ $(\text{\rm mod }P_2^* P_4^*)$ suitably later.
The letter $p$ will denote in the following an unspecified prime.
The choices of $v$ and $U$ in \eqref{eq:4.13} guarantee that
\beq
\label{eq:4.21}
(z \pm n, P_1^* P_3^*) = 1 \ \ \ (0 < n \leqslant U, \ n \neq 1)
\eeq
if and only if either
\beq
\label{eq:4.22}
n = pq^\alpha \prod_{i = 1}^m h_i^{\alpha_i} \ \ (\alpha \geqslant 0, \ \alpha_i \geqslant 0) \ \text{ and } \ n \geqslant R/2
\eeq
or
\beq
\label{eq:4.23}
n \ \text{ is composed only of primes } \ p \mid P_2^* q \prod_{i = 1}^m h_i .
\eeq

Our first (and main) goal in finding $z$ $\text{\rm mod } P^*$ will be to find residues $\alpha_p$ with $p \mid P_2^*$ so that choosing
\beq
\label{eq:4.24}
z \equiv \alpha_p \ (\text{\rm mod } p) \ \text{ for } \ p \mid P_2^*
\eeq
the condition $(z^* \pm h_i, P^*) = 1$ (cf.\ \eqref{eq:4.4}) should remain true and simultaneously we should have for as many as possible numbers $n$ of the form \eqref{eq:4.22} $(z + n, P^*) > 1$.
By a suitable choice of the parameter $w$ we can make the whole set of $n$'s satisfying \eqref{eq:4.23} relatively small, due to de Bruijn's result \cite{Bru} which we use in a weaker form,
proved by Rankin \cite{Ran1}.
The present form is Lemma~5 of \cite{Mai1}.

\begin{lemma}
\label{lem:1}
Let $\Psi(x,y)$ denote the number of positive integers $n \geqslant x$ which are composed only of primes $p \leqslant y$.
For $y \leqslant x$ and $y$ approaching to infinity with $x$, we have
\beq
\label{eq:4.25}
\Psi(x,y) \leqslant x \exp \left[ - \frac{\log_3 y}{\log y} \log x + \log_2 y + O \left(\frac{\log_2 y}{\log_3 y}\right)\right].
\eeq
\end{lemma}

Since $\Psi(x,y)$ is clearly monotonically increasing in $x$ we can estimate the number of $n$'s in \eqref{eq:4.23} from above by
\beq
\label{eq:4.26}
(\log U)^{m + 1} \Psi(U, w).
\eeq
Now, choosing
\beq
\label{eq:4.27}
w = \exp \bigl[\alpha(\log R \log_3 R / \log_2 R)\bigr]
\eeq
with a constant $\alpha$ to be determined later,
the quantity in \eqref{eq:4.26} is by \eqref{eq:4.13} clearly
\begin{align}
\label{eq:4.28}
&\ll U \exp \left[ - \frac{(1 + o(1)) \log_3 R \log R}{\alpha \log R \log_3 R/\log_2 R} + (1 + o(1)) \log_2 R + (m + 1) \log_2 R\right] \\
&\ll \frac{U}{\log^2 U} = o \left(\frac{R}{\log R}\right) \nonumber
\end{align}
if we chose, e.g.,
\beq
\label{eq:4.29}
\alpha = \frac1{m + 5} .
\eeq

This bound will be completely satisfactory for us.
So we have to concentrate on the numbers $n$ in \eqref{eq:4.22}.

The simplest strategy in choosing $\alpha_p$ $(\text{\rm mod }p)$ in \eqref{eq:4.24} would be to sieve out for the smallest $p > v$ the largest residue class $\text{\rm mod }p$ among the $N_0$ numbers of the form \eqref{eq:4.22} and repeat this consecutively for all primes in~$P_2^*$.
What makes our task more complicated, is the fact that we have to preserve the condition $(z + h_i, M) = 1$ in \eqref{eq:4.4}, so we have to choose for any $\widetilde p_j \in P_2^*$
\beq
\label{eq:4.30}
\alpha_{\widetilde p_j} \not\equiv -h_i \ (\text{\rm mod }\widetilde p_j) \ \text{ for } \ i = 1,2, \dots, m.
\eeq

First we observe that an easy calculation shows that the Prime Number Theorem implies
\begin{align}
\label{eq:4.31}
N_0 &= \sum_{\alpha, \alpha_1, \dots, \alpha_m \geqslant 0} \pi \left(\frac{U}{q^\alpha \prod\limits_{i = 1}^m {h_i}^{\alpha_i}}\right) - \pi \left(\frac{R}{2q^\alpha \prod\limits_{i = 1}^m {h_i}^{\alpha_i}}\right) \\
&\sim \frac{U}{\log U} \left(1 + \frac1q + \frac1{q^2} + \dots\right) \prod_{i = 1}^m \left(1 + \frac1{h_i} + \frac1{h_i^2} + \dots\right)
\nonumber\\
&\sim \frac{U}{\log U} \left(1 - \frac1q \right)^{-1} \prod_{i = 1}^m \left(1 - \frac1{h_i} \right)^{-1} \leqslant \frac{2U}{\log U}
\nonumber
\end{align}
if $C_3(m)$ in \eqref{eq:4.2} was chosen sufficiently large (taking into account also \eqref{eq:4.17}).

We will choose the residue classes $\alpha_p$ one by one for all primes from $v$ to $w$ and consider at the $i$th step the arising situation.
We will be left at the $j$th step with $N_j < N_0$ remaining values of $n$'s with \eqref{eq:4.22}.
If $N_j$ is at any stage of size $\leqslant \dfrac{R}{5 \log R}$, then we are ready.
Thus we can suppose $N_j > R/(5\log R)$.

Otherwise, if $\alpha$ and all $\alpha_i$ $(1 \leqslant i \leqslant m)$ are determined, then we have trivially in case of
\beq
\label{eq:4.32}
q^\alpha \prod_{i = 1}^m {h_i}^{\alpha_i} > \sqrt{U}
\eeq
at any rate altogether at most
\beq
\label{eq:4.33}
O\left(\sqrt{U}(\log U)^{m + 1}\right)
\eeq
numbers of the given form \eqref{eq:4.22}.

On the other hand, if \eqref{eq:4.32} is not true, then we have by the Brun--Titchmarsh theorem altogether at most
\beq
\label{eq:4.33masodik}
\frac{2m U}{\varphi\Bigl(\widetilde p_j q^\alpha \prod\limits_{i = 1}^m {h_i}^{\alpha_i}\Bigr) \log \sqrt{U}} \leqslant \frac{8U m}{\widetilde p_j q^{\alpha} \Bigl(\prod\limits_{i=1}^m {h_i}^{\alpha_i}\Bigr) \log U}
\eeq
numbers $n$ of form \eqref{eq:4.22} which lie in one of the $m$ bad residue classes $\{h_i\}_{i = 1}^m \text{ \rm mod }\widetilde p_j$.
Adding \eqref{eq:4.33masodik} up for all possible non-negative combination of $\alpha$, $\alpha_1, \dots, \alpha_m$ we get altogether (cf.\ \eqref{eq:4.31}) at most
\beq
\label{eq:4.34}
\frac{16m U}{\widetilde p_j \log U} < \frac{16 m c_0 (k_0) Rf(R)}{\widetilde p_j \log R} < \frac{N_j f(R)}{\widetilde p_j} < \frac{N_j}{4\log^2 R}
\eeq
bad $n$ values of the form \eqref{eq:4.22} by $\widetilde p_j > w$ and $N_j > R/(5 \log R)$.
This means that choosing from among the remaining $\widetilde p_j - m$ residue classes that one which sieves out the most elements from the remaining set of size $N_j$ we obtain for the size of the new remaining set
\begin{align}
\label{eq:4.36} 
N_{j + 1} &< N_j - \frac{N_j\left(1 - \frac1{4\log^2 R}\right)}{\widetilde p_j - m} \\
&< N_j \left(1 - \frac{1 - \frac1{4\log^2 R}}{\widetilde p_j}\right) \nonumber\\
&< N_j \left(1 - \frac1{\widetilde p_j}\right)^{1 - \log^{-2}R}.\nonumber
\end{align}
This means that by Mertens' theorem we obtain a final residual set of cardinality at most
\begin{align}
\label{eq:4.37}
N^* &< N_0 \prod_{\substack{v < p < w\\ p \neq q}} \left(1 - \frac1p\right)^{1 - \log^{-2}R} \sim N_0 \left(\frac{\log v}{\log w}\right)^{1 - \log^{-2}R} \\
&= N_0 \left(\frac{3 \log_2^2 R(1 + o(1))}{\alpha \log R \log_3 R}\right)^{1 - \log^{-2}R} \nonumber\\
&\leqslant 4 m N_0 \frac{\log_2^2 R}{\log R \log_3 R} \leqslant \frac{8m U}{\log U} \cdot \frac1{f(R)}\nonumber\\
&= \frac{8mc_0(k_0)R}{\log U} < \frac{8mc_0(k_0)R}{\log R} < \frac{\pi(R) - \pi(R/2)}{4} \nonumber
\end{align}
if the value of $c_0(k_0)$ in \eqref{eq:4.1} was chosen sufficiently small.

Formulas \eqref{eq:4.28} and \eqref{eq:4.31} mean that after the sieving with suitably chosen $\alpha_{\widetilde p_j}$ $(\text{\rm mod }\widetilde p_j)$, $\widetilde p_j \in P_2^*$ and taking into account that the set satisfying \eqref{eq:4.23} was at any rate small, all elements of \eqref{eq:4.22} and \eqref{eq:4.23} can be sieved out by one third of the primes from $P_4$ choosing for every remaining $n^* =
pq^\alpha \prod\limits_{i = 1}^m {h_i}^{\alpha_i}$ in \eqref{eq:4.22} and \eqref{eq:4.23} a separate $p^* \in P_4$, $p^* \neq p, q$ and $z \equiv - n^*$ $(\text{\rm mod }p^*)$, that is, $\alpha_{p*} \equiv - n^*$ $(\text{\rm mod }p^*)$.
This will assure $(z + n^*, M^*) > 1$.
With an other third of $p^* \mid P_4$ we can similarly sieve out the remaining negative $n^*$ values.
The only problem is that for any $p^* \in P_4$ used above, simultaneously with $z + n^* \equiv 0$ $(\text{\rm mod }p^*)$ we have to assure $z + h_i \not \equiv 0$ $(\text{\rm mod } p^*)$.
This is equivalent to that $n^* - h_i \equiv 0$ $(\text{\rm mod }p^*)$ is not allowed if we want to sieve out $z + n^*$ with $p^*$.
However, for every number $|n^*| \leqslant U$ and for every $i \in [1, m]$ the number $n^* - h_i$ has at most one prime divisor $> R/2 > \sqrt{U}$.
This means that for every $n^*$ we have at most $m$ forbidden primes to use to sieve out~$n^*$.
Thus, as the number of still ``abundant'' primes is at every step at least one third of all primes between $(R/2, R]$ we can for all (positive and negative) values of $n^*$ consecutively choose the primes $p^* \mid P_4^*$ nearly freely, just avoiding at most $m$ forbidden primes at each step, so that we would have, step by step for each $p^* \mid P_4^*$, $n^*$
\beq
\label{eq:4.38}
z + n^* \equiv 0 \ (\text{\rm mod } p^*), \ \ z + h_i \not \equiv 0 \ (\text{\rm mod } p^*) \ \ \ (i = 1,2,\dots, m).
\eeq

At the end of this procedure some $p^* \mid P_4$ will remain.
We can choose for these primes $\alpha_p^*$ freely with the only condition that
\beq
\label{eq:4.39}
\alpha_p^* \not\equiv - h_i \ (\text{\rm mod }p^*) \ \ \ (i = 1,2,\dots, m)
\eeq
should hold, in order to assure for the remaining primes $p^* \mid P_4$ also
\beq
\label{eq:4.40}
(z + h_i, p^*) = 1.
\eeq

Finally to determine $z$ uniquely $\text{\rm mod }M/q$ we choose for the primes $p = h_i$ \ $\alpha_p$ again essentially freely with the only condition
\beq
\label{eq:4.41}
\alpha_p \not \equiv - h_j \ (\text{\rm mod }p) \ \ \ (j = 1,2,\dots, m)
\eeq
which is again possible by $h_1 > m$.
In this way we will have finally for any $i = 1,2, \dots, m$ with $M^* = P^* = M/q = P/q$
\beq
\label{eq:4.42}
(z + h_i, P^*) = 1
\eeq
and this means that the relations \eqref{eq:4.4}--\eqref{eq:4.5} will be satisfied with $M^*$ in place of~$M$.

Now, we continue with the Maynard--Tao proof which we take over from \cite{May} with the additional changes executed in \cite{BFM}.
We will look for primes among the numbers
\beq
\label{eq:4.43}
\ell + h_1, \dots, \ell + h_m, \ \ \ell \equiv z \ (\text{\rm mod }W), \ \ W = M^* = P^*
\eeq
with the value $z$ found in the previous procedure, satisfying with the notation $\mathcal H_m = \{h_i\}_{i = 1}^m$
\begin{align}
\label{eq:4.44}
(z + h_i, W) &= 1 \ \ (i = 1,2,\dots, m)\\
\label{eq:4.45}
(z + \nu, W) &> 1 \ \text{ if } \ |\nu| \in \bigl[0, \log L f(\log L)\bigr] \setminus \mathcal H_m.
\end{align}
We will show that for $1 > L(k_0)$ we will find numbers
\beq
\label{eq:4.46}
\ell \in [L, 2L], \ \ \ell \equiv z\ (\text{\rm mod }W), \ \ \#\bigl\{\ell + h_i \in \mathcal P \ (1 \leqslant i \leqslant m)\bigr\} \geqslant k_0
\eeq
with \eqref{eq:4.43}--\eqref{eq:4.45} which will prove our Theorem~\ref{th:2}.

\begin{rema}
If there are additional primes among $\ell - h_i$ $(1 \leqslant i < m)$ this does not change anything since $h_i \leqslant C_4(m) \leqslant C_8(k_0)$.
The introduction of the new variables $\ell$ and $L$ is only necessary since we look for primes around $\ell$ instead of investigating $p_n$ and $d_n$, so essentially $\ell \sim n \log n$, $L \sim  N \log N$, $\log L \sim \log N$, so the size of gaps \eqref{eq:4.1} remains unchanged if we substitute $N$ by~$L$.
We remark that \eqref{eq:4.2} and \eqref{eq:4.17} assure
\beq
\label{eq:4.47}
p \,\bigg| \prod_{\substack{i, j = 1\\ i \neq j}}^m (h_i - h_j) \Longrightarrow p = {\mathcal O}_m(1), \ p \neq h_t \ (1 \leqslant t \leqslant m) \Longrightarrow p \nmid W,
\eeq
analogously to \eqref{eq:4.27}--\eqref{eq:4.31} of \cite{BFM}, which is actually the Maynard--Tao theorem.

We need also the modified Bombieri--Vinogradov theorem, Theorem~4.1 of \cite{BFM}, which is somewhat similar to, but stronger than Theorem~6 of \cite{GPY3}.
The greatest prime factor of $W$ is $\ll \log L$, so the level of smoothness of $W$ is much better than that $(L^\varepsilon)$ required by the condition of Theorem~10 of \cite{Cha}.
So the whole Theorem 4.1 of \cite{BFM} will remain true.
Further, we can leave out from the actual sieving procedure the possibly existing $q$ and its multiples.
This causes just a negligible change of size $\left(1 - \dfrac1q\right)^{O(1)} = 1 + O \left(\dfrac1{\log L}\right)$ in the weighted number of primes and in the sum of weights.
If we choose in \eqref{eq:4.3}
\beq
\label{eq:4.48}
R \leqslant c_9(k_0) \log L
\eeq
with a sufficiently small $c_9(k_0) > 0$, then we will have by the Prime Number Theorem
\beq
\label{eq:4.49}
W \leqslant L^{c_9(k_0)(1 + o(1))}
\eeq
so the whole Maynard--Tao procedure will remain valid, as in \cite{BFM}.
(The variable $R$ is here completely different from that in \cite{May} or \cite{BFM}.)
\end{rema}

Summarizing, we will obtain at least $k_0$ bounded prime gaps in intervals of type
\beq
\label{eq:4.50}
[\ell - h_m, \ell + h_m]
\eeq
and around them two intervals of size $(1 + o(1)) c_0(k_0)\log L \log_2 L \log_4 L/ (\log_3^2 L)$
containing only composite numbers (see \eqref{eq:4.7}) which prove Theorem~\ref{th:2}.

Finally, the fact that we left out the largest prime factor $q$ (and its multiples) of the possibly existing exceptional modulus $r$ yields an effective modified Bombieri--Vinogradov theorem and so an effective final Theorem~\ref{th:2}, consequently Theorems~\ref{th:1}, 1' and \ref{th:3} are effective too.

\bigskip
\noindent
{\bf Acknowledgement:} The author would like to express his sincere
gratitude to Imre Z. Ruzsa, who called his attention that a possible
combination of the methods of Erd\H{o}s--Rankin and Zhang--Maynard--Tao might
lead to stronger results about the ratio of consecutive primegaps than
those proved by the author in his earlier work \cite{Pin2}.

\noindent
J\'anos Pintz\\
R\'enyi Mathematical Institute\\
of the Hungarian Academy of Sciences\\
Budapest, Re\'altanoda u. 13--15\\
H-1053 Hungary\\
e-mail: pintz.janos@renyi.mta.hu

\end{document}